\documentclass{article}
\usepackage{amsmath,amscd,amssymb, latexsym, amsthm, multicol}
\usepackage[dvipdfmx]{graphicx,color}

\newtheorem{thm}{Theorem}
\newtheorem{prop}{Proposition}

\newtheorem{problem}{Problem}

\newtheorem{cor}{Corollary}

\theoremstyle{definition}

\def\Lower{\mathrm{Lower}}

\def\real{\mathbb{R}} 
\def\complex{\mathbb{C}} 
\def\field{\mathbb{K}}
\def\qed{\hfill $\Box$}

\begin{document}
\title{Lowerable vector fields\\ 
for a finitely ${\cal L}$-determined multigerm}
\author{Yusuke Mizota\footnote{Research Fellow DC2 of Japan Society for the Promotion of Science
\newline
} \ and Takashi Nishimura}
\date{}
\maketitle
\begin{abstract}  
We show that the module of lowerable vector fields for a finitely ${\cal L}$-determined 
multigerm is finitely generated in a constructive way. 
\end{abstract}

\section{Introduction}\label{intro}
Let $\field$ be $\real$ or $\complex$. 
Throughout this paper,  unless otherwise stated, all mappings are smooth 
(that is,  of class $C^{\infty}$ if $\field$ = $\real$ or holomorphic if $\field$ = $\complex$). 

Let $S$ be a subset consisting of $r$ distinct points in $\field^n$. 
A map-germ $f:(\field^n,   S)\rightarrow (\field^p,  0)$ is called a \textit{multigerm}. 
A multigerm $f:(\field^n,  S)\rightarrow (\field^p,  0)$ 
can be identified with $\{ f_k:(\field^n,  0)\rightarrow (\field^p,  0) |\, 1\leq k\leq r \}$.
Each $f_k$ is called a \textit{branch} of $f$. If $r=1$, $f$ is called a \textit{monogerm}. 

Let $C_{n,S}$ (resp., $C_{p,0}$) be the set of function-germs 
$(\field^n,  S)\rightarrow \field$ (resp., $(\field^p, 0)\rightarrow \field$),  
and let $m_{n,S}$ (resp., $m_{p,0}$) be the subset of $C_{n,S}$ (resp.,  $C_{p,0}$) consisting of function-germs 
$(\field^n,  S)\rightarrow (\field,  0)$ (resp.,  $(\field^p,  0)\rightarrow (\field,  0)$).  
For a non-negative integer $i$, let $m_{n,S}^i$ 
be the subset of $C_{n,S}$ consisting of function-germs 
$(\field^n,  S)\rightarrow \field$ such that the terms of the Taylor series of them up to ($i-1$) are zero.  
The sets $C_{n,S}$ and $C_{p,0}$ have natural $\field$-algebra structures. 
For a multigerm $f:(\field^n, S)\rightarrow (\field^p, 0)$, 
let $f^*:C_{p,0}\rightarrow C_{n,S}$ be the $\field$-algebra homomorphism defined by 
$f^*(\psi)=\psi \circ f$. 
Set $Q(f)=C_{n,S}/f^{*}m_{p,0}C_{n,S}$ and $\delta(f)= \dim_{\field} Q(f)$. 

For a map germ $f:(\field^n,  S)\rightarrow \field^p$,  
let $\theta_{n,S}(f)$ be the set of germs of vector fields along $f$. 
The set $\theta_{n,S}(f)$ has a natural $C_{n,S}$-module structure and is identified with the direct sum of $p$ copies of $C_{n,S}$.  
Put $\theta_{S}(n)=\theta_{n,S}(\mathrm{id}_{(\field^n,  S)})$ and 
$\theta_0(p)=\theta_{p, \{0\}}(\mathrm{id}_{(\field^p,  0)})$,  
where $\mathrm{id}_{(\field^n,  S)}$ (resp.,  $\mathrm{id}_{(\field^p,  0)}$) is the germ of the identity mapping of 
$(\field^n,  S)$ (resp.,  $(\field^p,  0)$).  
For a multigerm $f:(\field^n,  S)\rightarrow (\field^p,  0)$,  
following Mather \cite{ma}, define $tf$ : $\theta_{S}(n)\rightarrow \theta_{n,S}(f)$ (resp., $\omega f$ : $\theta_0(p)\rightarrow \theta_{n,S}(f)$) as 
$tf(\xi)=df\circ \xi$ (resp., $\omega f(\eta)=\eta\circ f$),
where $df$ is the differential of $f$. 
Following Wall \cite{wa}, set $T{\cal R}_e(f)=tf(\theta_{S}(n))$ and $T{\cal L}_e(f)= \omega f(\theta_0(p))$.  
For a multigerm $f:(\field^n,  S)\rightarrow (\field^p,  0)$, 
a vector field $\xi\in\theta_S(n)$ is said to be \textit{lowerable} if $df\circ\xi$ $\in T{\cal R}_e(f)\cap T{\cal L}_e(f)$. 
The set of lowerable vector fields is denoted by $\Lower(f)$.   
The set $\Lower(f)$ has a $C_{p,0}$-module structure via $f$.  
The notion of lowerable vector field, which was introduced by Arnol'd \cite{ar} for studying bifurcations of wave front singularities, is significant in Singularity Theory (for instance, see \cite{is}). 

In this paper, we investigate the following problem. 
\begin{problem}\label{prob1}
Let $f:(\field^n, S)\rightarrow (\field^p, 0)$ be a multigerm satisfying $\delta(f)<\infty$. Then, 
is $\Lower(f)$ finitely generated? In the case that $\Lower(f)$ is finitely generated, prove it in a constructive way. 
\end{problem}

Our first result is the following Proposition \ref{prop1}, which reduces Problem \ref{prob1} to the problem on $T{\cal R}_e(f)\cap T{\cal L}_e(f)$. 

\begin{prop}\label{prop1}
Let $f:(\field^n, S)\rightarrow (\field^p, 0)$ be a multigerm satisfying $\delta(f)<\infty$. Then, $tf$ is injective. 
\end{prop}
\noindent
In the complex analytic case, since $C_{p,0}$ is Noetherian and $T{\cal R}_e(f)\cap T{\cal L}_e(f)$
is a $C_{p,0}$-submodule of the finitely generated module  $\theta_{n,S}(f)$, it follows that $T{\cal R}_e(f)\cap T{\cal L}_e(f)$ is finitely generated. However, the algebraic argument gives 
no constructive proof. Moreover, in the real $C^{\infty}$ case, even finite generation of  $T{\cal R}_e(f)\cap T{\cal L}_e(f)$ seems to be open. 

The main purpose of this paper is to give a constructive proof of the following theorem, which works well in both the real $C^{\infty}$ case and the complex analytic case. 

\begin{thm}\label{thm1}
Let $f:(\field^n, S)\rightarrow (\field^p, 0)$ be a finitely ${\cal L}$-determined multigerm. 
Then, $T{\cal R}_e(f)\cap T{\cal L}_e(f)$ is finitely generated as a $C_{p,0}$-module via $f$.
\end{thm}
\noindent
Here, a multigerm $f:(\field^n, S)\rightarrow (\field^p, 0)$ is said to be \textit{finitely ${\cal L}$-determined}
if 
there exists a positive integer $\ell$  such that  
$m_{n,S}^\ell\theta_{n,S}(f) \subset T{\cal L}_e(f)$ holds. 
It is easily seen that if $f$ is finitely ${\cal L}$-determined, then $\delta(f)$ is finite. 
Thus, by combining Proposition \ref{prop1} and Theorem \ref{thm1}, we have the following partial affirmative answer to Problem \ref{prob1}. 

\begin{cor}\label{cor1}
Let $f:(\field^n, S)\rightarrow (\field^p, 0)$ be a finitely ${\cal L}$-determined multigerm. 
Then, $\Lower(f)$ is finitely generated as a $C_{p,0}$-module via $f$.
\end{cor}
\bigskip
In Section \ref{proof1} (resp., Section \ref{proof2}), Proposition \ref{prop1} (resp., Theorem  \ref{thm1}) is proved.  

\section*{Acknowledgement}
This work is partially supported by JSPS and CAPES under the Japan-Brazil research cooperative program. Y.M. is supported by Grant-in-Aid for JSPS Fellows Grant Number 251998.    

\section{Proof of Proposition \ref{prop1}}\label{proof1}
It is sufficient to show that if $f:(\field^n, 0)\rightarrow (\field^p, 0)$ is a 
monogerm satisfying $\delta(f)<\infty$, then $tf$ is injective. 
The property $\delta(f)<\infty$ implies that $f$ is finite to one (for instance, see \cite{GG}). 
Suppose that $tf (\xi) = 0$. Then it follows that $f$ is constant along any integral curve of $\xi$. Since $f$ is finite to one, it follows that any integral curve of $\xi$ must consist of only one point. This means that $\xi = 0$. Thus, $tf$ must be injective. This completes the proof. 
\qed

\section{Proof of Theorem \ref{thm1}}\label{proof2} 
Since $f$ is finitely ${\cal L}$-determined, 
there exists a positive integer $\ell$ such that 
\begin{equation}\label{eqn1}
m_{n,S}^\ell\theta_{n,S}(f) \subset T{\cal L}_e(f)
\end{equation}
and $\delta(f)<\infty$.  
Thus, we have that 
$Q(f_k)$ is a finite dimensional vector space for any $k$ ($1\leq k \leq r$). Set $\delta(f_k) = m_k$ and  $Q(f_k)=\{[\phi_{k,1}],\ldots,[\phi_{k,m_k}]\}_\field$, where $\phi_{k,j}\in C_{n,\{0\}}$ and 
$[\phi_{k,j}]=\phi_{k,j}+f_{k}^{*}m_{p, 0}C_{n,\{0\}}$ for any $j$ ($1\leq j \leq m_k$). 
 
Now, we find a finite set of generators for $T{\cal R}_e(f)\cap T{\cal L}_e(f)$. 
Let $(x_1,\ldots, x_n)$ (resp.,  $(X_1, \ldots, X_p)$) be the standard local coordinates of $\field^n$ (resp.,  $\field^p$) at the origin. 
Then, by the preparation theorem, $T{\cal R}_e(f_k)$ may be expressed as follows: 
\[
\left.
\left\{
\left(\begin{array}{ccc}
\displaystyle{\frac{\partial(f_{k,1})}{\partial x_1}}&\cdots&\displaystyle{\frac{\partial(f_{k,1})}{\partial x_n}} \\
\vdots&&\vdots \\
\displaystyle{\frac{\partial(f_{k,p})}{\partial x_1}}&\cdots&\displaystyle{\frac{\partial(f_{k,p})}{\partial x_n}}
\end{array}\right)
\left(\begin{array}{c}
\displaystyle{\sum_{1\leq j\leq m_k}(\psi_{k,1,j}\circ f_k)\phi_{k,j}} \\
\vdots \\
\displaystyle{\sum_{1\leq j\leq m_k}(\psi_{k,n,j}\circ f_k)\phi_{k,j}}
\end{array}\right)
\right|\,
\psi_{k,i,j}\in C_{p,0}
\right\}
,
\]
where $f_k=(f_{k,1},\ldots,f_{k,p})$.
It can be simplified as follows:  
\[
\left.
\left\{
\sum_{1\leq i\leq n, 1\leq j\leq m_k }^{}(\psi_{k,i,j}\circ f_k)\xi_{k,i,j} 
\right|\, 
\psi_{k,i,j}\in C_{p,0} 
\right\}, 
\]
where $\xi_{k,i,j}$ is the transposed matrix of 
\[
\left(\frac{\partial (f_{k,1})}{\partial x_i}\phi_{k,j},\ldots,\frac{\partial (f_{k,p})}{\partial x_i}\phi_{k,j} \right). 
\]
Note that $\xi_{k,i,j} \in T{\cal R}_e(f_k)$. For simplicity, we abbreviate $\displaystyle{\sum_{1\leq i\leq n, 1\leq j\leq m_k }}$ to $\displaystyle{\sum_{i, j}}$.
Let $\bar{\eta}=(\bar{\eta_1},\ldots,\bar{\eta_r})$ be an element of 
 $T{\cal R}_e(f)\cap T{\cal L}_e(f)$. 
Then, $\bar{\eta}$ may be expressed as follows: 
\[
\bar{\eta}=(\bar{\eta}_1,\ldots,\bar{\eta}_r)=\left(\sum_{i,j}(\psi_{1,i,j}\circ f_1)\xi_{1,i,j},\ldots,\sum_{i,j}(\psi_{r,i,j}\circ f_r)\xi_{r,i,j}\right).
\]
For a $p$-tuple of non-negative integers $\alpha=(\alpha_1,\ldots,\alpha_p)$, 
set $|\alpha|=\alpha_1+\cdots+\alpha_p$, $f_k^{\alpha}=f_{k,1}^{\alpha_1}\cdots f_{k,p}^{\alpha_p}$ and $X^{\alpha}=X_1^{\alpha_1}X_2^{\alpha_2}\cdots X_p^{\alpha_p}$.
The function-germ $\psi_{k,i,j} \in C_{p,0}$ can be written as the sum of a polynomial of total degree less than or equal to ($\ell-1$) and an element in $m^\ell_{p,0}$ as follows: 
\[
\psi_{k,i,j}(X_1,\ldots,X_p)=\sum_{0 \leq |\alpha| \leq \ell-1}c_{k,i,j,\alpha}X^{\alpha}
+\sum_{|\alpha|=\ell}\tilde{\psi}_{k,i,j,\alpha}X^{\alpha},
\]
where $c_{k,i,j,\alpha} \in \field$ and $\tilde{\psi}_{k,i,j,\alpha} \in C_{p,0}$. 
Note that $\ell$ is the positive integer given in $(\ref{eqn1})$. 
Then, we have the following: 
\[
\bar{\eta}_k=\sum_{i,j} \, \sum_{0 \leq |\alpha| \leq \ell-1}c_{k,i,j,\alpha}(f_k^{\alpha}\xi_{k,i,j})
+\sum_{i,j} \, \sum_{|\alpha|=\ell}\tilde{\psi}_{k,i,j,\alpha}(f_k^{\alpha}\xi_{k,i,j}). 
\]
Set $\bar{\xi}_{k,i,j,\alpha}= (0,\dots,0,\underbrace{f_k^{\alpha}\xi_{k,i,j}}_{k\mbox{-th component}},0,\ldots,0)$.
Note that $\bar{\xi}_{k,i,j,\alpha}\in T{\cal R}_e(f)$. 
Then, we have the following:
\[
\bar{\eta} = \sum_{1\leq k\leq r}\sum_{i,j}\sum_{0 \leq |\alpha| \leq \ell-1}c_{k,i,j,\alpha}\bar{\xi}_{k,i,j,\alpha} 
+\sum_{1\leq k\leq r}\sum_{i,j}\sum_{|\alpha| = \ell}(\tilde{\psi}_{k,i,j,\alpha} \circ f)\bar{\xi}_{k,i,j,\alpha}. 
\]
We define finite sets $L$ and $H$ as follows:
\[
L=\left\{\bar{\xi}_{k,i,j,\alpha} \left|\, 
 0 \leq |\alpha| \leq \ell-1, 1\leq k \leq r, 1\leq i \leq n, 1\leq j \leq m_k\right\}\right.,
\]
\[
H=\left\{\bar{\xi}_{k,i,j,\alpha} \left|\, 
|\alpha|=\ell, 1\leq k \leq r, 1\leq i \leq n, 1\leq j \leq m_k\right\}\right..
\]
Then, $H \subset T{\cal R}_e(f)\cap T{\cal L}_e(f)$ by (\ref{eqn1}). 
Therefore, 
\[
\sum_{1\leq k\leq r}\sum_{i,j}\sum_{0 \leq |\alpha| \leq l-1}c_{k,i,j,\alpha}\bar{\xi}_{k,i,j,\alpha}
\]
belongs to $V=T{\cal R}_e(f)\cap T{\cal L}_e(f) \cap L_{\field}$.  
The set $V$ is a finite dimensional vector space. 
Set $\dim_{\field} V = m$. Then, there exist $\underline{\xi}_1,\ldots,\underline{\xi}_m \in T{\cal R}_e(f)\cap T{\cal L}_e(f)$ such that $V=\langle\underline{\xi}_1,\ldots,\underline{\xi}_m\rangle_{\field}$.
It is clear that $V \subset \langle\underline{\xi}_1,\ldots,\underline{\xi}_m\rangle_{f^*C_{p,0}}$.
Therefore, 
\[
\bar{\eta} \in \langle\underline{\xi}_1,\ldots,\underline{\xi}_m\rangle_{f^*C_{p,0}}+H_{f^*C_{p,0}}. 
\]
Hence, 
\[
T{\cal R}_e(f)\cap T{\cal L}_e(f)
\subset 
\langle\underline{\xi}_1,\ldots,\underline{\xi}_m\rangle_{f^*C_{p,0}}+H_{f^*C_{p,0}}. 
\]
Since $\{\underline{\xi}_1,\ldots,\underline{\xi}_m\} \cup H$ is contained in $T{\cal R}_e(f)\cap T{\cal L}_e(f)$, 
the converse inclusion also holds.  
Therefore, $T{\cal R}_e(f)\cap T{\cal L}_e(f)$ is finitely generated as a $C_{p,0}$-module via $f$.
This completes the proof. 
\qed

Graduate School of Mathematics, Kyushu University, 744, Motooka, Nishi-ku, Fukuoka 819-0395, JAPAN. \newline e-mail: \texttt{y-mizota@math.kyushu-u.ac.jp}
\newline

Research Group of Mathematical Sciences, Research Institute of Environment and Information Sciences, Yokohama National University, Yokohama 240-8501, JAPAN. \newline e-mail: \texttt{nishimura-takashi-yx@ynu.jp}


\begin{thebibliography}{99}
\bibitem{ar} V.  I.  Arnol'd,  \textit{Wave front evolution and equivariant Morse lemma},  Commun Pure Appl.  Math.,  \textbf{29} (1976),
557--582. 

\bibitem{GG} M. Golubitsky and V. Guillemin, \textit{Stable Mappings and Their Singularities}, Graduate Texts in Math., \textbf{14}, Springer-Verlag, New York (1973).

\bibitem{is} G. Ishikawa,  \textit{Openings of differentiable map-germs and unfoldings}, Topics on Real and Complex Singularities, Proceedings of the 4th Japanese-Australian Workshop (JARCS4), Kobe 2011, World Scientific (2014), 87--113. 

\bibitem{ma} J.  Mather,  \textit{Stability of $C^\infty$ mappings,  I\!I\!I.  Finitely determined map-germs},  Publ.  Math.  Inst.  Hautes Etudes Sci.,  \textbf{35} (1969), 127--156. 

\bibitem{wa} C.  T.  C.  Wall,  \textit{Finite determinacy of smooth map-germs},  Bull.  London Math.  Soc.,  \textbf{13} (1981), 481--539. 

\end{thebibliography}
\end{document}